\documentclass[12pt]{amsart}
\allowdisplaybreaks
\usepackage{amssymb}
\usepackage{graphicx}
\usepackage{fullpage}
\usepackage[pdfencoding=auto,colorlinks,linkcolor=,citecolor=]{hyperref}
\usepackage{cite}
\usepackage{dsfont}
\usepackage{eucal}
\usepackage[all,cmtip]{xy}

\newcommand{\Coloc}{\mathsf{Coloc}}
\newcommand{\DD}{\mathsf{D}}
\newcommand{\Db}{\DD^{\mathsf{b}}}

\newcommand{\Dsg}{\DD_{\mathsf{sg}}}

\newcommand{\iEnd}{\mathcal{E}\mathit{nd}}

\newcommand{\Ext}{\mathsf{Ext}}
\newcommand{\tExt}{\mathsf{E\widehat{\phantom{\dot{}}x}t}}

\newcommand{\Hom}{\mathsf{Hom}}
\newcommand{\iHom}{\mathcal{H}\mathit{om}}

\newcommand{\ik}{\mathsf{ik}}
\newcommand{\ind}{\mathop{\rm ind}}

\newcommand{\kG}{\kk G}
\newcommand{\kk}{k}
\newcommand{\KInj}{\mathsf{KInj}}
\newcommand{\KInjb}{\KInj^{\mathsf{b}}}
\newcommand{\KacInj}{\mathsf{K^{ac}Inj}}
\newcommand{\Loc}{\mathsf{Loc}}
\newcommand{\Locb}{\Loc^{\mathsf b}}

\newcommand{\Mod}{\mathsf{Mod}}
\newcommand{\pk}{\mathsf{pk}}
\newcommand{\phat}{{}^{^\wedge}_p}

\newcommand{\Perf}{\mathsf{Perf}}

\newcommand{\res}{\mathop{\rm res}}
\newcommand{\sHom}{\underline{\Hom}}
\newcommand{\Soc}{\mathsf{Soc}}
\newcommand{\Spec}{\mathsf{Spec}\,}
\newcommand{\stmod}{\mathsf{stmod}}
\newcommand{\StMod}{\mathsf{StMod}}
\newcommand{\Thick}{\mathsf{Thick}}
\newcommand{\tk}{\mathsf{tk}}
\newcommand{\tM}{\mathsf{t}M}

\newcommand{\bF}{\mathbb F}

\newcommand{\bZ}{\mathbb Z}

\newcommand{\cC}{\mathcal C}

\newcommand{\cF}{\mathcal F}

\newcommand{\cT}{\mathcal{T}}

\newcommand{\fm}{\mathfrak{m}}

\renewcommand{\le}{\leqslant}
\renewcommand{\ge}{\geqslant}

\makeindex
\numberwithin{equation}{section}

\theoremstyle{plain}
\newtheorem{lemma}[equation]{Lemma}
\newtheorem{theorem}[equation]{Theorem}
\newtheorem{proposition}[equation]{Proposition}
\newtheorem{corollary}[equation]{Corollary}

\theoremstyle{definition}
\newtheorem{conjecture}[equation]{Conjecture}

\newtheorem{definition}[equation]{Definition}

\theoremstyle{remark} 
\newtheorem{remark}[equation]{Remark}

\author{David J. Benson} 
\address{Institute of Mathematics \\ 
University of Aberdeen \\ 
Aberdeen AB24 3UE \\ 
United Kingdom}
\email{d.j.benson@abdn.ac.uk}

\author{John Greenlees} 
\address{Mathematics Institute, Zeeman Building, University of
  Warwick, Coventry CV4 7AL, United Kingdom}
\email{John.Greenlees@warwick.ac.uk}

\title{Modules with finitely generated cohomology, and singularities of $C^*BG$}

\begin{document}

\begin{abstract}
Let $G$ be a finite group and $\kk$ a field of characteristic $p$. We
conjecture that if $M$ is a $\kG$-module with $H^*(G,M)$ finitely
generated as a module over $H^*(G,\kk)$ then 
as an element of the
stable module category $\StMod(\kG)$, $M$ is contained in the thick
subcategory generated by the finitely generated $\kG$-modules and the
modules $M'$ with $H^*(G,M')=0$.

We show that this is equivalent to a conjecture of the second author
about generation of the bounded
derived category of cochains $C^*(BG;\kk)$, and
we prove the conjecture in the case where the centraliser of every
element of $G$ of order $p$ is $p$-nilpotent. In this case some
stronger statements are true, that probably fail for more general
finite groups.
\end{abstract}

\maketitle

\section{Introduction}

Let $G$ be a finite group and $\kk$ a field of characteristic $p$. We
are interested in some conjectures about objects with finitely
generated cohomology in three closely related situations: the stable
module category $\StMod(\kG)$, the homotopy category of complexes of
injective $\kG$-modules $\KInj(\kG)$, and the derived category
$\DD(C^*BG)$. Throughout this paper, we write
$C^*BG$ for $C^*(BG;\kk)$, the cochains on $BG$ with coefficients in
$\kk$.

In each of these situations, we write down a stronger conjecture, which
sometimes holds but is probably too strong in most cases, and a weaker
conjecture which we hope is always true. We prove that the stronger
conjectures for a given $\kG$ in the three contexts are all equivalent, and that the
weaker conjectures for $\kG$ in the three contexts are also all
equivalent.  

We also show that the weaker conjectures are equivalent
to a conjecture of the second author. This is the
statement that $\Db(C^*BG)$ is generated by the cochains $C^*BS$,
where $S$ is a Sylow $p$-subgroup of $G$.
It is not hard to deduce by formal arguments that the thick
subcategory generated by $C^*(BS)$ contains $C^*(BH)$ for all 
subgroups $H$ of $G$ 
(see also Conjecture~\ref{conj:weak-DC}). Accordingly the conjecture
is really about showing that all finitely generated modules over
$C^*(BG)$ are finitely built from the most obvious ones.

We prove the stronger conjectures in the case where the centraliser of
every element of order $p$ in $G$ is $p$-nilpotent, and a separate
paper~\cite{Benson/Carlson:bc19} proves a few more cases of the 
weaker conjectures, but the general
case of the weaker conjectures remains open.

If $M$ and $N$ are $\kG$-modules, we say that a homomorphism $M \to N$ is
a \emph{cohomology equivalence} if the induced map on Tate cohomology
is an isomorphism. This is not an equivalence relation, but
we complete it to one, and then talk of
$M$ and $N$ being cohomologically
equivalent. We say that $M$ is a \emph{no cohomology module} 
if its Tate cohomology vanishes.

\begin{conjecture}\label{conj:strong-stmod}
If $\hat H^{\ge 0}(G,M)$ is finitely generated as a module over $H^*(G,\kk)$ then
$M$ is cohomologically equivalent to a finitely generated $\kG$-module.
\end{conjecture}

This is the stronger conjecture for the stable module category.
We write $\StMod(\kG)$ for the stable 
category of $\kG$-modules, a compactly generated triangulated category 
whose compact objects $\stmod(\kG)$ consist of the thick subcategory 
of modules stably isomorphic to finitely generated modules. We write $\StMod^0(\kG)$ for the thick
subcategory of no cohomology modules, and $\StMod^b(\kG)$ for the thick
subcategory of $\kG$-modules $M$ such that $\hat H^{\ge 0}(G,M)$ is finitely
generated as a module over
$H^*(G,\kk)$. Conjecture~\ref{conj:strong-stmod}
states that every object in $\StMod^b(\kG)$ is equivalent in the
Verdier quotient $\StMod^b(\kG)/\StMod^0(\kG)$ to the image of an
object in $\stmod(\kG)$.

We are able to prove this conjecture under some rather restrictive
hypotheses, but in general it is
not at all clear that the images of objects in $\stmod(\kG)$ even form
a triangulated subcategory of the Verdier quotient. So
the following weaker conjecture for the
stable module category states just that their images generate, and seems more
likely to be true in general.

\begin{conjecture}\label{conj:weak-stmod}
$\StMod^b(\kG)$ is equal to the thick subcategory 
\[ \Thick(\stmod(kG),\StMod^0(\kG)) \] 
of $\StMod(kG)$ generated by the finitely
generated modules and the no cohomology modules.
\end{conjecture}

Next we discuss the context of $\KInj(\kG)$.
We write 
$\KInj(\kG)$ for the homotopy category of complexes of injective
$\kG$-modules, studied in~\cite{Benson/Krause:2008a}. 
The full subcategory $\KacInj(\kG)$ of acyclic complexes
is equivalent to $\StMod(\kG)$ via Tate resolutions. The thick
subcategory of compact objects in $\KInj(\kG)$ is equivalent to
$\Db(\kG)$ via semi-injective resolutions, and we identify it with
$\Db(\kG)$.

We write
$\KInj^0(\kG)$ for the no cohomology objects in $\KInj(\kG)$, namely the  
localising subcategory of objects $X$ with $H^*(G,X)=0$. We write
$\KInjb(\kG)$ for the thick subcategory of objects such that
$H^*(G,X)$ is finitely generated as a module over $H^*(G,\kk)$. The
following are the analogues of Conjectures~\ref{conj:strong-stmod} 
and~\ref{conj:weak-stmod}, and we prove that they are equivalent to
the corresponding statements about $\StMod(\kG)$.

\begin{conjecture}\label{conj:strong-Db}
If $X$ is in $\KInjb(\kG)$ then $X$ is cohomologically equivalent to
an object in $\Db(\kG)$.
\end{conjecture}

This is the same as saying that the image of $X$ in the Verdier
quotient $\KInjb(\kG)/\KInj^0(\kG)$ is the image of an object in
$\Db(\kG)$. Again, it is not even clear that these images form a
triangulated subcategory, so the following is the weaker conjecture in
this context.

\begin{conjecture}\label{conj:weak-Db}
$\KInjb(\kG)$ is equal to the thick subcategory 
\[ \Thick(\Db(\kG),\KInj^0(\kG)) \] 
of $\KInj(\kG)$ generated
by the bounded derived category and the no cohomology objects.
\end{conjecture}

Next,
let $\DD(C^*BG)$ be the derived category of cochains on $BG$.
We can view this in several ways. If we think of $C^*BG$ as a
differential graded algebra, then we take the homotopy category of
differential graded modules, and we invert the quasi-isomorphisms.
If we think of $C^*BG$ as an $A_\infty$ algebra~\cite[\S4]{Keller:2001a},
then we take the homotopy category of $A_\infty$ modules. There is
then no need to invert quasi-isomorphisms, as they are automatically
invertible. If we think of $C^*BG$ as an
algebra over the sphere spectrum~\cite[Chapter~III]{Elmendorf/Kriz/Mandell/May:1997a},
then we take the category
of modules, and invert the weak equivalences. These three points of
view give rise to equivalent triangulated categories.

Write $\Db(C^*BG)$ for the thick subcategory of objects $X$ in
$\DD(C^*BG)$ such that
\[ H^*(X)=\Hom_{\DD(C^*BG)}(C^*BG,X) \] 
is finitely generated as a module over $H^*BG=H^*(BG;\kk)$. This
definition is justified and discussed in Greenlees and
Stevenson~\cite{Greenlees/Stevenson:2020a}. The following are the
stronger and weaker conjecture in this context.

\begin{conjecture}\label{conj:strong-DC}
If $Y$ is in $\Db(C^*BG)$ then there is an object $X$ in $\Db(\kG)$
such that $Y$ is equivalent to $C^*(BG;X)$.
\end{conjecture}

\begin{conjecture}\label{conj:weak-DC}
$\Db(C^*BG)$ is generated as a thick subcategory of $\DD(C^*BG)$ by the objects
$C^*(BG;X)$ for $X$ in $\Db(\kG)$.
\end{conjecture}

We show that Conjecture~\ref{conj:weak-DC} is equivalent to
the following conjecture of the second author.

\begin{conjecture}\label{conj:C*BS}
The bounded derived category $\Db(C^*BG)$ is generated by the
module $C^*BS$, where $S$ is a Sylow $p$-subgroup of $G$.
\end{conjecture}

\begin{remark}
A minor variant of Conjectures~\ref{conj:weak-DC} and~\ref{conj:C*BS}
is that $\Db(C^*BG)$ 
is generated by the objects $C^*(EG\times_G Y)$ for finite
$G$-complexes $Y$.
\end{remark}

Here is a diagram of the implications between the various conjectures.
\[ \xymatrix{\ref{conj:strong-stmod} \ar@{<=>}[r]\ar@{=>}[d] &
\ref{conj:strong-Db} \ar@{<=>}[r]\ar@{=>}[d] & 
\ref{conj:strong-DC}\ar@{=>}[d]\\
\ref{conj:weak-stmod}\ar@{<=>}[r] & 
\ref{conj:weak-Db}\ar@{<=>}[r] &
\ref{conj:weak-DC}\ar@{<=>}[r] & 
\ref{conj:C*BS}} \]
The vertical implications are obvious. The left hand horizontal equivalences are
proved in Theorem~\ref{th:equivalences}, while the middle ones
follow from Propopsition~\ref{pr:DbC*BG}, and the right hand one is
proved in Theorem~\ref{th:C*BS-equivalent}.

At the end of
Section~\ref{se:nucleus}, we prove the  following.

\begin{theorem}\label{th:C*BS-nonucleus}
  Conjectures~\ref{conj:strong-stmod}--\ref{conj:C*BS} all
  hold for $G$ if the centraliser of every element
  of order $p$ in $G$ is $p$-nilpotent.
\end{theorem}

In a separate paper~\cite{Benson/Carlson:bc19},
it is proved that the weaker conjectures in the bottom row of the above
diagram are true for some cases not covered by Theorem~\ref{th:C*BS-nonucleus}.\bigskip

Finally, we explain the role of the nucleus, a concept introduced
in~\cite{Benson:1995a,Benson/Carlson/Robinson:1990a}, and we
conjecture that this is the support of the singularity category of
$C^*BG$. We prove that if Conjecture~\ref{conj:C*BS} holds for a group
$G$ then the support of $\Dsg(C^*BG)$ is contained in the nucleus,
but equality seems more difficult to prove except in the cases where
the nucleus is either $\varnothing$ or $\{0\}$, see Theorem~\ref{th:triv-nucleus}.\bigskip

\noindent
{\bf Notation.} In a tensor triangulated category $\cT$ with shift $\Sigma$, we write
$\Hom(X,Y)$ for the internal Hom functor, adjoint to tensor
product, whose existence is guaranteed by Brown representability.
We write $\Hom^n(X,Y)$ for $\Hom(X,\Sigma^n Y)$, so that $\Hom^*(X,Y)$
is the graded $\Hom$ space. We write $\Hom_\cT(X,Y)$ for the
morphisms in the category $\cT$.

If objects in $\cT$ are chain complexes, we write $\iHom^*_{\cT}(X,Y)$
for the complex of homomorphisms, and note that its degree zero
cohomology is $\Hom_{\cT}(X,Y)$. Then  $\iEnd^*_{\cT}(X)=
\iHom^*_{\cT}(X,X)$ is a differential
graded algebra, and $\iHom^*_{\cT}(X,Y)$ is an
$\iEnd^*_{\cT}(Y)$--$\iEnd^*_{\cT}(X)$ bimodule. This makes $\cT$ the
homotopy category of a differential graded category.\bigskip

\noindent
{\bf Acknowledgement.} We would like to thank the Hausdorff Institute
of Mathematics in Bonn for its hospitality and funding during the
programme HIM-Spectral-22, while much of this
research was carried out. We would also like to thank Jon Carlson and
Ran Levi for conversations relevant to this work.

\section{Modules with finitely generated cohomology}

\begin{definition}
Recall that a $\kG$-module is projective if and only if it is injective.
We write $\StMod(\kG)$ for the stable category of $\kG$-modules,
and $\stmod(\kG)$ for the full subcategory of finitely generated
$\kG$-modules. These are triangulated categories in which the graded
homomorphisms from $M$ to $N$ are the Tate groups $\tExt^*_{\kG}(M,N)$. Recall also
that $\hat H^i(G,M)$ denotes $\tExt^i_{\kG}(\kk,M)$.
\end{definition}

Making the definition at the beginning of the introduction more
precise, we have the following.

\begin{definition}
If
$M$ and $N$ are $\kG$-modules, we say that a homomorphism $M \to N$ is
a \emph{cohomology equivalence} if the induced map $\hat H^*(G,M) \to
\hat H^*(G,N)$ is an isomorphism, or equivalently if
$\hat H^{\ge 0}(G,M) \to \hat H^{\ge 0}(G,N)$ is an isomorphism, see
Proposition~\ref{pr:gaps} and Corollary~\ref{co:gaps} below. This is
not an equivalence relation, but we complete it to one.
Note that we do not require this to hold
for subgroups of $G$. We say that $M$ and $N$ are cohomologically
equivalent if they are equivalent with respect to this equivalence
relation. In other words, there is a zig-zag of cohomological
equivalences relating $M$ and $N$.

We say that $M$ is a \emph{no cohomology module}
if the Tate cohomology $\hat H^*(G,M)=0$, or equivalently if
$\hat H^{\ge 0}(G,M)=0$.
\end{definition}

From the point of view of this paper, finite $p$-groups are not so
interesting, because the following theorem shows that
Conjecture~\ref{conj:strong-stmod} holds for these.

\begin{theorem}\label{th:p-group-stmod}
Suppose that $G$ is a finite $p$-group. If $M$ is a $\kG$-module such
that $H^0(G,M)$ is finite dimensional
then $M$ is finitely generated. If $\hat H^n(G,M)$ is finite dimensional
for some $n\in\bZ$
then $M$ is a direct sum of a finite dimensional module and a
projective module, and hence equivalent to a finite dimensional module
in the stable module category.
\end{theorem}
\begin{proof}
Recall that for a finite $p$-group in characteristic $p$, there is
only one simple module, namely the trivial $\kG$-module $\kk$.
Its projective cover and injective hull are both $\kG$, and every
injective is a direct sum of copies of $\kG$.

Let $I$ be the injective hull of $M$. Then the inclusion $M \to I$
induces an isomorphism $H^0(G,M) \to H^0(G,I)$. But $I$ is a direct
sum of copies of $\kG$, so the dimension of $I$ is $|G|.\dim H^0(G,I)$.
So if $H^0(G,M)$ is finite dimensional then $I$ is finite dimensional, and hence so is $M$.

Similarly, if $\hat H^n(G,M)$ is finite dimensional with $n\in\bZ$ then in 
a minimal Tate resolution $\cdots \to I_{-1} \to I_0 \to I_1\to \cdots$ of $M$, 
the dimension of $I_n$ is $|G|.\dim\hat H^n(G,M)$, and so the
dimension of $\Omega^{-n}M$ is finite. Finally, $M$ is a direct sum of
$\Omega^n\Omega^{-n}M$ and a projective module.
\end{proof}

The following is useful in studying no cohomology modules.

\begin{definition}
  Let $n_p(G)$ be the minimum over all homogeneous systems of parameters
$\zeta_1,\dots,\zeta_r$ for $H^*(G,\bF_p)$ of the expression
$1+\sum_{i=1}^r(|\zeta_i|-1)$. This is a positive integer depending
only on $p$ and $G$.
\end{definition}
  
\begin{proposition}\label{pr:gaps}\ 
\begin{enumerate} 
\item
If $\hat H^i(G,M)$ vanishes for $n_p(G)$ successive values of $i$ then
it vanishes for all $i$.
\item
If $\hat H^i(G,M)$ is finite dimensional for $n_p(G)$ successive values
of $i$ then it is finite dimensional for all $i$.
\end{enumerate}
\end{proposition}
\begin{proof}
A spectral sequence argument described in Theorem~2.4 of
Benson and Carlson~\cite{Benson/Carlson/Robinson:1990a} shows that the
first part holds, and an obvious 
modification of the same argument shows that the second part holds.
\end{proof}

\begin{corollary}\label{co:gaps}
If $f\colon M\to N$ is a $\kG$-module homomorphism inducing
isomorphisms $\hat H^i(G,M) \to \hat H^i(G,N)$ for $n_p(G)+1$
consecutive values of $i$ then it $f$ is a cohomology equivalence.
\end{corollary}
\begin{proof}
Complete to a triangle $M \to N \to M' \to \Omega^{-1}M$ in the stable
module category, look at the long exact sequence in cohomology, 
and apply part (1) of Proposition~\ref{pr:gaps} to the module $M'$.
\end{proof}

\begin{remark}
Recall that $H^n(G,M)$ agrees with $\hat H^n(G,M)$ for $n>0$. In
degree zero, we have a surjective map $H^0(G,M) \to \hat H^0(G,M)$
whose kernel is the ``projective'' fixed points of $M$, namely the
image of $\sum_{g\in G}g \in \kG$ on $M$. Expressing $M$ as a direct
sum of a projective module $P$ and a module $M'$ with no projective
summands, the kernel is the fixed points of $G$ on $P$ and the image
is the fixed points of $G$ on $M'$.

Combining this with Proposition~\ref{pr:gaps}, we see that if
$H^*(G,M)=0$ then $\hat H^*(G,M)=0$, and that if $\hat H^*(G,M)=0$
then $H^*(G,M)$ is concentrated in degree zero. Similarly, if
$H^*(G,M)$ is finitely generated as a module over $H^*(G,\kk)$
then $\hat H^{\ge 0}(G,M)$ is
also finitely generated, while if $\hat H^{\ge 0}(G,M)$ is finitely
generated then $M$ is a direct sum of a projective module $P$ and a
module $M'$ such that $H^*(G,M')$ is finitely generated.
So we may as well work
in the stable module category, which means that
we should ignore projective summands.
\end{remark}

Recall that Tate duality in $\StMod(\kG)$ says the following.

\begin{proposition}\label{pr:Tate-StMod}
 If $N$ is a finitely
generated $\kG$-module then for an arbitrary $\kG$-module $M$ we have
\[ \sHom^{-n-1}_{\kG}(M,N) \cong \Hom_\kk(\sHom^n_{\kG}(N,M),\kk). \]
\end{proposition}

The Brown--Comenetz dual $M^*$ of a 
module $M$ in the stable module category 
is defined by the isomorphism
\[ \Hom_\kk(\sHom_{\kG}(\kk,-\otimes M),\kk) \cong
  \sHom_{\kG}(-,M^*). \]
It follows from Tate duality that the Brown--Comenetz dual of a 
$\kG$-module $M$ in $\StMod(\kG)$ is 
\[ M^*\cong \Omega \Hom_k(M,k). \]
The Spanier--Whitehead dual $M^\vee$, on the other hand, is the function object
$\Hom_\kk(M,\kk)$ characterised by the adjunction
\[ \sHom_{\kG}(- \otimes M,\kk) \cong \sHom_{\kG}(-,\Hom_\kk(M,\kk)). \]
Here, the action of $\kG$ on $\Hom_k(M,k)$ is given via the 
antipode $\sigma$ of $\kG$ taking each group element to its inverse. 
See for example Section~8 of \cite{Benson/Iyengar/Krause:2011b}.\medskip

No cohomology modules are both the left and the right perp of the trivial module, in the
sense of the following definition and lemma.

\begin{definition}
If $\cC$ is an object or a collection of objects in a triangulated category $\cT$, we
write $\cC^\perp$ for the subcategory consisting of objects
with no non-zero maps from $\cC$, and ${}^\perp\cC$ for the localising
subcategory consisting of objects with no non-zero maps to $\cC$.
If $\cT$ has products then $\cC^\perp$ is colocalising, while if $\cC$
has coproducts then ${}^\perp\cC$ is localising.
\end{definition}

\begin{lemma}
If $N$ is a finitely generated $\kG$-module then $N^\perp$ and
${}^\perp N$ are equal as subcategories of $\StMod(\kG)$, and are both
localising and colocalising.
\end{lemma}
\begin{proof}
This follows directly from Tate duality.
\end{proof}

As an example, we have $\kk^\perp={}^\perp\kk=\StMod^0(\kG)$.
We have a recollement of triangulated categories
\[ \xymatrix{\kk^\perp \ar[r]|(0.4)u & \StMod(\kG) 
\ar[r]|(0.45)v\ar@<1ex>[l]^(0.6){u_\lambda}\ar@<-1ex>[l]_(0.6){u_\rho} & 
\StMod(\kG)/\kk^\perp \ar@<1ex>[l]^(0.55){v_\lambda}
\ar@<-1ex>[l]_(0.55){v_\rho},} \]
where the right hand side is the Verdier quotient.

Recall that $\StMod^b(\kG)$ is the subcategory of $\StMod(\kG)$
consisting of modules $M$ such that $\hat H^{\ge 0}(G,M)$ is finitely generated
as a module over $H^*(G,\kk)$. 

\begin{lemma}
$\StMod^b(\kG)$ is a thick subcategory of $\StMod(\kG)$.
\end{lemma}
\begin{proof}
$\StMod^b(\kG)$ is clearly closed under direct summands.
It follows from part (2) of Proposition~\ref{pr:gaps} that $M$ is in
$\StMod(kG)$ if and only if $\Omega(M)$ is, and then the long exact
sequence in cohomology shows that it is a triangulated subategory.
\end{proof}

Next we study duality for $\StMod^b(\kG)$.

\begin{proposition}\label{pr:dual}
If $M$ is in $\StMod^b(\kG)$ then so is the dual $M^\vee$.
\end{proposition}
\begin{proof}
We examine the local cohomology spectral
sequence~\cite{Greenlees:1995a}
\[ E_2^{s,t}=(H^{s}_\fm H^*(G,M))^t \Rightarrow H_{-s-t}(G,M). \]
By definition of Tate cohomology, the right hand side is
$\hat H^{s+t-1}(G,M)$ when $s+t\le 0$. 

Now suppose that $H^*(G,M)$ is finitely generated as a module over
$H^*(G,\kk)$. Then the $E_2$ term is Artinian as a module over
$H^*(G,\kk)$, and hence so are $H_*(G,M)$ and $\hat H^{\le 0}(G,M)$. In particular, it is finite
dimensional in each degree. So
by Tate duality, this is a shift of the graded dual 
of $H^*(G,M^\vee)$ as $H^*(G,\kk)$-modules. So the latter is
Noetherian as a graded module over $H^*(G,\kk)$.
\end{proof}

\begin{lemma}\label{le:dual}
Suppose that a $\kG$-module $M$ has no projective summands. 
Then neither does the dual $M^\vee=\Hom_k(M,k)$.
\end{lemma}
\begin{proof}
The problem we have to contend with here is that the kernel of the
action of $\kG$ on the dual $M^\vee$ is not the same as the kernel of
the action on $M$, because of the intervention of the antipode
$\sigma\colon\kG\to\kG$, $\sigma(g)=g^{-1}$, in the definition of the
dual. One way to get around this is to consider the double dual
$M^{\vee\vee}\cong M^{**}$. Since $\sigma^2$ is the identity, the
natural inclusion $M \to M^{\vee\vee}$ does not involve the antipode, and
the annihilators of $M$ and $M^{\vee\vee}$ are equal.

Since $kG$ is self-injective, a module $M$ has no projective summands
if and only if the socle of the group algebra $\Soc(\kG)$ acts as zero
on $M$. If this holds then $\Soc(\kG)$ also acts as zero on
$M^{\vee\vee}$, so $M^{\vee\vee}$ has no projective summands. On the
other hand, the dual of a projective module is projective, so if
$M^\vee$ has a projective summand then so does $M^{\vee\vee}$.
\end{proof}

\begin{proposition}\label{pr:double-dual}
For $M$ in $\StMod^b(\kG)$, the natural map to the double dual 
\[ M \to M^{**} \cong M^{\vee\vee} \] 
is a cohomology isomorphism.
\end{proposition}
\begin{proof}
We may suppose without loss of generality that $\Omega^{-i}M$ has no projective
summands. Then by Lemma~\ref{le:dual}, nor does $(\Omega^{-i}M)^*$ or
$(\Omega^{-i}M)^{**}=(\Omega^iM^*)^*$, so we may take
\[ \Omega^{-i}(M^{**})=(\Omega^iM^*)^*=(\Omega^{-i}M)^{**}. \] 
So 
\[ \Hom_{\kG}(\kk,\Omega^{-i}(M^{**})) \cong (\Hom_{\kG}(\Omega^iM^*,k))^*
\cong (\Hom_{\kG}(k,\Omega^{-i}M))^{**}. \]
Then since there are no projective summands, we have 
$\Hom_{\kG}(k,\Omega^{-i}M) \cong \hat H^i(G,M)$. In particular,
if $\hat H^i(G,M)$ is finite dimensional then the map 
\[ \hat H^i(G,M)\to \hat H^i(G,M^{**})\cong \hat H^i(G,M)^{**} \] 
is an isomorphism as required.
\end{proof}

\begin{corollary}
If there is a counterexample to Conjecture~\ref{conj:strong-stmod} then there is a pure
injective counterexample.
\end{corollary}
\begin{proof}
This follows from Proposition~\ref{pr:double-dual}, 
since a module is pure injective if
and only if it is isomorphic to a direct summand of the dual of a module.
\end{proof}

\section{\texorpdfstring{The category $\KInj(\kG)$}
{The category KInj(kG)}}

\begin{definition}
Let $\KInj(\kG)$ be the homotopy category of complexes of injective
 $\kG$-modules. This is a tensor triangulated category. Its tensor
identity $\ik$ is an injective resolution of $\kk$ as a $\kG$-module.
We also write $\pk$ for a projective resolution of $\kk$ and $\tk$ for
a Tate resolution of $\kk$. These sit in a triangle
\[ \pk \to \ik \to \tk \to \Sigma\pk. \]
Note that we index objects in $\KInj(\kG)$ cohomologically, so that
$\ik$ lives in positive degrees while $\pk$ lives in negative
degrees. As a result, the endomorphism rings of $\ik$ and $\pk$, which
are both the cohomology ring, live in positive
degree. 
If $X$ is an object in $\KInj(\kG)$, we write 
\[ H^*(G,X) =
  \Hom^*_{\KInj(\kG)}(\ik,X)\cong H^*(\Hom^*_{\kG}(\kk,X)). \]
\end{definition}

See Benson and Krause~\cite{Benson/Krause:2008a} for background on
this category. Its compact objects are the semi-injective resolutions 
of finite complexes of finitely generated $\kG$-modules, and these form
a thick subcategory equivalent to the bounded derived category
$\Db(\kG)$.

Every complex $X$ of injective $\kG$-modules is a direct sum of a split
exact sequence of injectives and a complex $Y$ that is \emph{minimal} in
the sense that the socle in each degree is in the kernel of the
differential. The inclusion $Y \to X$ is a homotopy equivalence,
so objects in $\KInj(\kG)$ may be assumed to be minimal when
desired. For a minimal complex $X$ we have $H^n(G,X)\cong
\Hom_{\kG}(\kk,X^n)$. 

We have an embedding of $\StMod(\kG)$ into $\KInj(\kG)$, taking
a module $M$ to a Tate resolution $\tM$ of $M$. This sits in a
recollement
\[ \xymatrix{\StMod(\kG)\simeq\KacInj(\kG)\ar[r] &
\KInj(\kG) \ar@<1ex>[l]\ar@<-1ex>[l]\ar[r] &
\DD(\kG)\ar@<1ex>[l]\ar@<-1ex>[l]} \]
where $\KacInj(\kG)$ is the full subcategory of acyclic complexes, 
and  $\DD(\kG)$ is the derived category. The left adjoints in this
recollement preserve compactness and give us functors on compact objects
\[ \stmod(\kG) \leftarrow \Db(\kG) \leftarrow \Perf(\kG) \]
where $\Perf(\kG)$ consists of the perfect complexes. This gives
Rickard's equivalence
\[ \stmod(\kG)\simeq\Db(\kG)/\Perf(\kG), \]
see Theorem~2.1 of~\cite{Rickard:1989a}.

Given any object $X$ in $\KInj(\kG)$, the left adjoints give a
triangle
\[ \pk \otimes X \to X \to \tk \otimes X \to \Sigma\pk\otimes X \]
while the right adjoints give a triangle
\[ \Hom(\tk,X) \to X \to \Hom(\pk,X) \to \Sigma\Hom(\tk,X). \]

Tate duality in $\KInj(\kG)$ says the following; see for example
Krause and Le~\cite{Krause/Le:2006a}.

\begin{proposition}\label{pr:Tate-KInj}
If $C$ is in $\Db(\kG)\subseteq\KInj(\kG)$ then for all $X$ in
$\KInj(\kG)$ we have
\[ \Hom_{\kG}(X,\pk\otimes C) \cong
  \Hom_\kk(\Hom_{\kG}(C,X),\kk). \]
\end{proposition}

It follows from Tate duality that 
the Brown--Comenetz dual is given by  
\[ X^*= \Hom(X,\pk) \cong \Hom_{\kk}(X,\kk) \] 
where $\pk$ is a projective resolution of $\kk$,
while the Spanier--Whitehead dual is given
by
\[ X^\vee = \Hom(X,\ik). \]

\begin{remark}
If we take the vector space dual of a Tate resolution of a
$\kG$-module $M$, we get  a Tate resolution of the Brown--Comenetz
dual $M^*=\Omega\Hom_{\kk}(M,\kk)$, which is a shift of a Tate resolution of the
Spanier--Whitehead dual $M^\vee=\Hom_\kk(M,\kk)$. This explains 
the shift in numbering between Tate duality in $\StMod(\kG)$ and in
$\KInj(\kG)$ (compare Proposition~\ref{pr:Tate-StMod}
with Proposition~\ref{pr:Tate-KInj}).
\end{remark}

\begin{lemma}\label{le:kG/S-gen}
The category $\Db(\kG)$ is generated by the simple $\kG$-modules. It
is also generated by the permutation module $\kk(G/S)$ on the cosets
of a Sylow $p$-subgroup $S$.
\end{lemma}
\begin{proof}
Generation by the simple modules follows from the fact that every
finitely generated module has a finite filtration with simple filtered
quotients. In particular, for a
$p$-group $S$, $\Db(\kk S)$ is generated by $\kk$.

For a  general finite group, given an object $X$ in $\Db(\kG)$, 
the restriction $\res_{G,S}(X)$ is generated by $\kk$, and so inducing
back up to $G$, the module
$\ind_{S,G}\res_{G,S}(X)$ is generated by $\ind_{S,G}(\kk)=\kk(G/S)$.
The second statement now follows from the fact that $X$ is isomorphic 
to a direct summand of $\ind_{S,G}\res_{G,S}(X)$.
\end{proof}

\section{\texorpdfstring{The conjectures in $\KInj(\kG)$}
{The conjectures in KInj(kG)}}

We have $\ik^\perp={}^\perp\ik=\KInj^0(\kG)$, the
no cohomology objects in $\KInj(\kG)$.
We have a recollement of triangulated categories
\[ \xymatrix{\KInj^0(\kG) \ar[r]|u & \KInj(\kG) 
\ar[r]|(0.37)v\ar@<1ex>[l]^{u_\lambda}\ar@<-1ex>[l]_{u_\rho} & 
\KInj(\kG)/\KInj^0(\kG) \ar@<1ex>[l]^(0.63){v_\lambda}
\ar@<-1ex>[l]_(0.63){v_\rho}} \]
where the right hand side is the Verdier quotient. The image of the
right adjoint
$v_\rho$ is the localising subcategory $\Loc(\ik)$ of $\KInj(\kG)$
generated by $\ik$, while the image of the left adjoint $v_\lambda$ is the colocalising
subcategory $\Coloc(\ik)$ cogenerated by $\ik$. 

Given an object $X$
in $\KInj(\kG)$, this gives us a localisation triangle
\[ v_\rho v X \to X \to u_\rho u X \to \Sigma v_\rho v X \]
where $v_\rho v X$ is in $\Loc(\ik)$, $v_\rho v X \to X$ is a cohomology equivalence and $u_\rho u X$
is in $\KInj^0(\kG)$. Similarly, it gives us a colocalisation triangle
\[ u_\lambda u X \to X \to v_\lambda v X \to \Sigma u_\lambda u X \]
where $v_\lambda v X$ is in $\Coloc(\ik)$, $X \to v_\lambda v X$ is a cohomology equivalence and 
$u_\lambda u X$ is in $\KInj^0(\kG)$.

\begin{definition}
We define $\KInjb(\kG)$ to be the thick subcategory of
$\KInj(\kG)$ consisting of objects $X$ whose cohomology
$H^*(G,X)$ is finitely generated as a module over $H^*(G,k)$.
Thus $\Db(\kG)\subseteq\KInjb(\kG)$ and $\KInj^0(\kG)\subseteq\KInjb(\kG)$.
We define $\Locb(\ik)$ to be the thick subcategory 
$\Loc(\ik)\cap \KInjb(\kG)\subseteq\KInj(\kG)$. So $v_\rho$ induces an
equivalence
\begin{equation}\label{eq:KInjb/KInj0}
\KInjb(\kG)/\KInj^0(\kG)\xrightarrow{\sim}\Locb(\ik).
\end{equation}
\end{definition}

We now examine Conjectures~\ref{conj:strong-Db} and~\ref{conj:weak-Db}.
Again, finite $p$-groups are not so interesting from the point of view
of this paper, because the following theorem shows that
Conjecture~\ref{conj:strong-Db} holds in this case.

\begin{theorem}\label{th:Kinjb-pgroup}
Suppose that $G$ is a finite $p$-group. Then $\KInjb(\kG)=\Db(\kG)$.
\end{theorem}
\begin{proof}
  Suppose that $X$ is an object in $\KInjb(\kG)$.
Thinking of $H^*(G,X)$ as $\Hom_{\KInj(\kG)}(\ik,X)$, it follows from
the hypothesis that if $Y$ is in $\Thick(\ik)$ then
$\Hom_{\KInj(\kG)}(Y,X)$ is finitely generated as a module over
$H^*(G,\kk)$. In the case of a finite $p$-group, the module $\kG$,
regarded as an element of $\Db(\kG)$ concentrated in a single degree,
is in $\Thick(\ik)$, as it is built from $|G|$ copies of $\kk$. So 
$\Hom_{\KInj(\kG)}(\kG,X)$ is finitely generated. But it is also
killed by the ideal of positive (cohomological) degree elements, so it is finite
dimensional. But $\Hom_{\KInj(\kG)}(\kG,X)$ is the cohomology of $X$ as
a complex. So $X$ is exact in all but finitely many degrees, and is
hence equivalent to a semi-injective resolution of an object in $\Db(\kG)$.
\end{proof}

\begin{proposition}\label{pr:whittle-down}
Given a complex $X$ in $\KInjb(\kG)$ there is a subcomplex $Y$
such that the inclusion $Y\to X$  is a cohomology equivalence,
and such that $Y$ is zero in all sufficiently large negative (cohomological)
degrees, and finite dimensional in all degrees.
\end{proposition}
\begin{proof}
Let $X$ be an object in $\KInjb(\kG)$.
Since $H^*(G,X)$ is finitely generated over $H^*(G,\kk)$, 
there exists $n\in \bZ$ such that
$H^i(G,X)=0$ for $i< n$. 
Assuming without loss of generality that $X$ is minimal, we have
$\Hom_{\kG}(\kk,X^i)=0$ for $i< n$. Then the (brutal)
truncation $X^{< n}$ is in $\KInj^0(\kG)$, and the inclusion 
$X^{\ge  n}\to X$ is a cohomology equivalence. 

We may now assume that $X^i=0$ for $i<n$. Then there is a finite
dimensional summand $Y^n$ of $X^n$ such that the inclusion
induces an isomorphism $\Hom_{\kG}(\kk,Y^n) \to \Hom_{\kG}(\kk,X^n)$. 
Next, there is a finite dimensional summand $Y^{n+1}$ of $X^{n+1}$
which contains the image of $Y^n$ under the differential, and such
that $\Hom_{\kG}(\kk,Y^{n+1}) \to \Hom_{\kG}(\kk,X^{n+1})$ is an
isomorphism. Continuing this way, we construct a subcomplex $Y$ of $X$
which is a summand in each degree, so that it is a complex of
injectives, and such that the quotient $X/Y$ is in
$\KInj^0(\kG)$. Then the inclusion $Y \to X$ is a cohomology equivalence.
\end{proof}

\begin{proposition}\label{pr:pk-Artinian}
Given a complex $X$ in $\KInjb(\kG)$ 
we have 
$H^n(G,\pk \otimes X)=0$ for $n$ sufficiently large positive.
\end{proposition}
\begin{proof}
By Proposition~\ref{pr:whittle-down}, we may assume 
that $X$ is zero in all sufficiently
negative degrees. So we have a local cohomology spectral sequence
\[ (H^s_\fm H^*(G,X))^t \Rightarrow H_{-s-t}(G,X), \]
see Greenlees~\cite[\S2]{Greenlees:1995a}. Here $H_*(G,X)$ denotes the
homology of the $G$-fixed points on $\pk\otimes X$, so that
$H_{-n}(G,X) \cong H^n(G,\pk\otimes X)$. Since $H^*(G,X)$ is
finitely generated over $H^*(G,\kk)$, each $H^s_\fm H^*(G,X)$ is
Artinian. The spectral sequence has a finite number of columns
(between the depth and the Krull dimension of $H^*(G,X)$) and so
$H^*(G,\pk\otimes X)$ is Artinian, and zero in all large enough
positive degrees.
\end{proof}

\begin{theorem}\label{th:tX}
Suppose that $X$ is in $\KInjb(\kG)$.
Then the map $\ik\to\tk$ induces an isomorphism $H^n(G,X)\to
H^n(G,\tk\otimes X)$ for all large enough $n$.
\end{theorem}
\begin{proof}
By
Proposition~\ref{pr:pk-Artinian}, $H^n(G,\pk \otimes X)=0$ for $n$
sufficiently large positive. The long exact sequence of the
triangle 
\[ \pk \otimes X \to X \to \tk\otimes X \to \Sigma\pk\otimes X \]
now proves the theorem.
\end{proof}

\begin{corollary}\label{co:tX}
We have the following.
\begin{enumerate}
\item
Suppose that $X$ is in $\KInj^0(\kG)$. If we regard $\tk\otimes X$
as a Tate resolution of a $\kG$-module $M$, then $M$ is in
$\StMod^0(\kG)$.
\item
Suppose that $X$ is in $\KInjb(\kG)$.
If we regard $\tk\otimes X$ as a Tate resolution of a $\kG$-module
$M$, then $M$ is in $\StMod^b(\kG)$.
\end{enumerate}
\end{corollary}
\begin{proof}
Statement (1) follows from Theorem~\ref{th:tX} and
Proposition~\ref{pr:gaps}\,(1), while statement
(2) follows from Theorem~\ref{th:tX} and
Proposition~\ref{pr:gaps}\,(2).
\end{proof}

It is clear that Conjecture~\ref{conj:strong-stmod} implies
Conjecture~\ref{conj:weak-stmod} and that
Conjecture~\ref{conj:strong-Db} implies Conjecture~\ref{conj:weak-Db}.
We now prove the following.

\begin{theorem}\label{th:equivalences}
Conjecture~\ref{conj:strong-stmod} is equivalent to
Conjecture~\ref{conj:strong-Db}, and Conjecture~\ref{conj:weak-stmod}
is equivalent to Conjecture~\ref{conj:weak-Db}.
\end{theorem}
\begin{proof}
We begin by showing that Conjecture~\ref{conj:strong-stmod} implies
Conjecture~\ref{conj:strong-Db} and Conjecture~\ref{conj:weak-stmod}
implies Conjecture~\ref{conj:weak-Db}. 
  
Let $X$ be an object in $\KInjb(\kG)$.
Replacing $X$ by the
cohomologically equivalent $v_\rho v X$, which we may do in both
Conjecture~\ref{conj:strong-Db} and Conjecture~\ref{conj:weak-Db},
we may suppose that $X$ is in $\Loc(\ik)$. Using
Proposition~\ref{pr:whittle-down} and
shifting if necessary we may assume that $X^n=0$ for $n<0$.  
By Corollary~\ref{co:tX}, $\tk \otimes X$ is a Tate resolution of a
$\kG$-module $M$ in $\StMod^b(\kG)\cap \Loc(\ik)$. 

If
Conjecture~\ref{conj:strong-stmod} holds, there is a finitely generated
$\kG$-module $M'\in \stmod(\kG)$ such that $M\to M'$ is a cohomology
equivalence; whereas if Conjecture~\ref{conj:weak-stmod} holds then
$M$ is a module in $\Thick(\stmod(\kG),\StMod^0(\kG))$, and we set
$M'=M$.

The module $M'$ is
only well defined up to adding and removing projective summands, but we
regard it as an object in $\Mod(\kG)$.
Then the complex $\ik \otimes M'$ is an injective resolution of $M'$, and 
may be regarded as an object in 
$\Db(\kG)\subseteq\KInj(\kG)$, respectively $\Thick(\Db(\kG),\KInj^0(\kG))$,
giving us a diagram
\[ \xymatrix@=6mm{\pk \otimes X \ar[r] & X \ar[r]\ar@{.>}[d]
 & \tk\otimes X = \tM\ar[d] \\
&\ik \otimes M' \ar[r]& \tk \otimes M' \ar[r] & \Sigma \pk \otimes M'.} \]
Since $X^n=0$ for $n<0$, and $(\Sigma\pk\otimes M')^n=0$ for $n\ge 0$, 
 there are no no-zero maps from $X$ to $\Sigma
\pk \otimes M'$. So we can fill
in the dotted arrow to make the diagram above commutative. This dotted
arrow is an isomorphism in cohomology in sufficiently large degrees.
Complete this to a triangle
\[ Y \to X \to \ik\otimes M' \to \Sigma Y. \]
Then $Y$ has cohomology only in finitely many degrees, and is in
$\Loc(\ik)$, so it is in $\Db(\kG)$. Since $\ik \otimes M'$ is  in
$\Db(\kG)$, respectively $\Thick(\Db(\kG),\KInj^0(\kG))$ 
it follows that $X$ is. This completes the proof that
Conjecture~\ref{conj:strong-stmod} implies
Conjecture~\ref{conj:strong-Db} and that
Conjecture~\ref{conj:weak-stmod} imples Conjecture~\ref{conj:weak-Db}.

Conversely, 
let $M$ be a module in $\StMod^b(\kG)$. We may suppose that $M$ has no
projective summands. Then $\ik \otimes M$ is an
object in $\KInjb(\kG)$. If Conjecture~\ref{conj:strong-Db} holds,
then 
there is an object $X$ in $\Db(\kG)$
that is cohomology equivalent to $\ik\otimes M$; whereas if
Conjecture~\ref{conj:weak-Db} holds, then $\ik\otimes M$ is in
$\Thick(\Db(\kG),\KInj^0(\kG))$, and we set $X=\ik\otimes M$.
Then $\tk\otimes X$
is in $\stmod(\kG)$, respectively $\Thick(\stmod(\kG),\StMod^0(\kG))$, 
and is cohomology equivalent to $M$.
\end{proof}

\section{\texorpdfstring{The category $\DD(C^*BG)$}{The category D(C*BG)}}

An injective resolution of the trivial module 
$\ik$ is the tensor identity in $\KInj(\kG)$, and we
have a quasi-isomorphism $\iEnd_{\kG}(\ik)\simeq C^*BG$ given by the
Rothenberg--Steenrod construction, where
$\iEnd_{\kG}(\ik)$ is the DG algebra of endomorphisms of $\ik$,
see for example~\cite[\S4]{Benson/Krause:2008a}.
If $X$ is an object in $\KInj(\kG)$, we write 
$\cF(X)=\iHom_{\kG}(\ik,X)$ for the differential graded Hom object from $\ik$ to $X$ in
$\KInj(\kG)$. 
Thus $\iHom_{\kG}(\ik,X)$ can be regarded as an object in
$\DD(\iEnd_{\kG}(\ik))$, and corresponds to the object $C^*(BG;X)$ in 
$\DD(C^*BG)$.

\begin{theorem}\label{th:DC*BG}
The functor 
\[ C^*(BG;-) \colon \KInj(\kG) \to \DD(C^*BG) \]
is essentially surjective, and induces equivalences of triangulated
categories  
\[ \Loc(\ik)\simeq\KInj(\kG)/\KInj^0(\kG)\simeq \DD(C^*BG). \] 
\end{theorem}
\begin{proof}
See for example \cite{Benson/Krause:2008a}.
\end{proof}

Since $\KInj^0(\kG)$ is contained in $\KInjb(\kG)$, 
the image of $\KInjb(\kG)$ consists of the $C^*BG$-modules with
finitely generated cohomology. Greenlees and
Stevenson~\cite{Greenlees/Stevenson:2020a} give good justifications for
the following definition.

\begin{definition}
The \emph{bounded derived category} $\Db(C^*BG)$ is defined to be the
full thick subcategory of $\DD(C^*BG)$ consisting of objects $X$ such
that $H^*(BG;X)$ is finitely generated as a module over $H^*(BG;\kk)$.
\end{definition}

\begin{proposition}\label{pr:DbC*BG}
The functor $\KInj(\kG) \to \DD(C^*BG)$ sending $X$ to $C^*(BG;X)$
induces an equivalence of triangulated categories
\[ \KInjb(\kG)/\KInj^0(\kG) \to \Db(C^*BG). \]
\end{proposition}
\begin{proof}
This follows from Theorem~\ref{th:DC*BG}.
\end{proof}

\begin{theorem}\label{th:DbC*BG2}
\begin{enumerate}
\item[\rm (i)]
Conjecture~\ref{conj:strong-Db} is equivalent
to Conjecture~\ref{conj:strong-DC}, and 
Conjecture~\ref{conj:weak-Db} is equivalent to
Conjecture~\ref{conj:weak-DC}.
\item[\rm (ii)]
If Conjecture~\ref{conj:strong-DC} holds then $\Db(C^*BG)$ is
equivalent to the
Verdier quotient of $\Db(\kG)$ by the thick subcategory of objects with
no cohomology.
\end{enumerate}
\end{theorem}
\begin{proof}
Both parts follow immediately from Proposition~\ref{pr:DbC*BG}.
\end{proof}

Now let $\ik(G/S)$ be the image in $\KInj(\kG)$ of the permutation
module $\kk(G/S)$ in $\Db(\kG)$. This is an injective resolution of
$\kk(G/S)$, or equivalently the induced complex $\ind_{S,G}(\ik)$.

\begin{lemma}\label{le:kG/S}
We have $C^*(BG;\ik(G/S))\cong C^*BS$.
\end{lemma}
\begin{proof}
Under the equivalence $\DD(\iEnd_{\kG}(\ik))\simeq\DD(C^*BG)$, the
object 
\[ \iHom_{\kG}(\ik,\ik(G/S)) \cong \iEnd_{\kk S}(\ik) \] 
goes to $C^*BS$ as a $C^*BG$-module.
\end{proof}

\begin{theorem}\label{th:C*BS-equivalent}
Conjectures~\ref{conj:C*BS} and \ref{conj:weak-DC} are equivalent.
\end{theorem}
\begin{proof}
By Lemma~\ref{le:kG/S}, we have $C^*(BG;\ik(G/S))\cong C^*BS$, so 
since $\ik(G/S)$ is in $\Db(\kG)$,
Conjecture~\ref{conj:C*BS} implies
Conjecture~\ref{conj:weak-DC}. Conversely, by Lemma~\ref{le:kG/S-gen},
$\ik(G/S)$ generates $\Db(\kG)$, so if the objects $C^*(BG;X)$ with
$X$ in $\Db(\kG)$ generate $\Db(C^*BG)$ then so does the object
$C^*(BG;\ik(G/S))=C^*BS$. So Conjecture~\ref{conj:weak-DC} implies
Conjecture~\ref{conj:C*BS}. 
\end{proof}

\section{The nucleus}\label{se:nucleus}

We use a slight variation of the definition of \emph{nucleus} from the
papers of Benson, Carlson and Robinson~\cite{Benson/Carlson/Robinson:1990a} 
and Benson~\cite{Benson:1995a}
to accommodate the fact that we are working with the bounded derived
category $\Db(\kG)$ rather than the stable module category 
$\stmod(\kG)$. So we define the
nucleus for $G$ to be the subset $\Theta_G$ of $V_G=\Spec H^*BG$ with the
following characterisations.
\begin{enumerate}
\item $\Theta_G$ is the union of the supports $V_G(Y)$ of those
  objects $Y$ in $\Db(B_0(\kG))$ such that $H^*(G,Y)=0$. Here
  $B_0(\kG)$ is the principal block of $\kG$.
\item $\Theta_G$ is the union of the images of $\res^*_{G,E}\colon V_E
  \to V_G$ as $E$ ranges over those elementary abelian $p$-subgroups
of $G$  such that $C_G(E)$ is not $p$-nilpotent.
\end{enumerate}
The equivalence of these characterisations is conjectured
in~\cite{Benson/Carlson/Robinson:1990a} and proved in
Benson~\cite{Benson:1995a}. 

\begin{remark}
If $G$ itself is $p$-nilpotent this definition gives
$\Theta_G=\varnothing$, while if $G$ is not $p$-nilpotent but the
centraliser of every element of order $p$ in $G$ is $p$-nilpotent then $\Theta_G=\{0\}$.
This distinction was not necessary in the context of
$\stmod(B_0(\kG))$ in the papers cited,
but is necessary in the current context of $\Db(B_0(\kG))$.
\end{remark}

\begin{theorem}\label{th:nucleus}
The category $\Db(B_0(\kG))$ is generated by $\kk$ and 
finitely generated modules $M$ with $V_G(M)\subseteq \Theta_G$.
\end{theorem}
\begin{proof}
This is proved in Theorem~1.2 of~\cite{Benson:1995a}.
\end{proof}

The following is essentially contained in~\cite{Benson:1995a}.

\begin{theorem}\label{th:zero-nucleus}
If $\Theta_G\subseteq\{0\}$, i.e., if the
centraliser of every element of order $p$ in $G$ is $p$-nilpotent,
every non-zero module in $\StMod(B_0(\kG))$ has non-zero cohomology.
The category $\StMod(\kG)$
decomposes as a direct sum of triangulated categories 
$\Loc(\kk) \oplus \StMod^0(\kG)$.
\end{theorem}
\begin{proof}
If $\Theta_G=\varnothing$ or $\Theta_G=\{0\}$ then by
Theorem~\ref{th:nucleus}, $\Db(B_0(\kG))$ is generated by $\kk$ and
projective modules in $B_0(\kG)$, so $\stmod(B_0(\kG))=\Thick(\kk)$
and $\StMod(B_0(\kG))=\Loc(\kk)$.
Then $\StMod^0(\kG)$ is the sum of the other blocks.
\end{proof}

\begin{corollary}\label{co:zero-nucleus}
If $\Theta_G=\varnothing$ or $\Theta_G=\{0\}$ then
Conjectures~\ref{conj:strong-stmod}--\ref{conj:C*BS} hold for $\kG$.
\end{corollary}
\begin{proof}
Since we've shown that Conjecture~\ref{conj:strong-stmod} implies the
remaining conjectures, we prove Conjecture~\ref{conj:strong-stmod}
in this case.
If $\Theta_G\subseteq\{0\}$ and $M$ is a module in $\StMod^f(\kG)$ then
by Theorem~\ref{th:zero-nucleus}, $M$
decomposes as a direct sum of a module in $\Loc(\kk)$ and
a no cohomology module. So $M$ is cohomologically equivalent to a 
module in $\Loc(\kk)$, and hence without loss of generality $M$ is
in $\Loc(\kk)$. So $\Ext^*_{\kG}(C,M)$ is finitely generated
over $H^*(G,\kk)$ whenever $C$ is in the thick subcategory generated
by $\kk$ and the no cohomology modules.

In particular, if $S$ is a Sylow $p$-subgroup of $G$ then the
permutation module $\kk(G/S)$ decomposes as a direct sum
of a module in $\Thick(\kk)$ and a no cohomology module.
So $\Ext^*_{\kG}(\kk(G/S),M)$ is finitely generated over $H^*(G,\kk)$.

By Shapiro's lemma $\Ext^*_{\kG}(\kk(G/S),M)\cong H^*(S,M)$ is finitely
  generated as a module over $H^*(G,\kk)$ acting through restriction
  to $P$, and hence also as a module over $H^*(S,\kk)$. So by
Theorem~\ref{th:p-group-stmod}, the restriction of $M$ to $S$ is finite dimensional plus
projective, and cohomologically equivalent to a finite dimensional
module. But $M$ is isomorphic to a direct summand of $M$ restricted to
$S$ induced back up to $G$, and hence cohomologically equivalent to a
finite dimensional module.
\end{proof}

\section{\texorpdfstring{The singularity category of $C^*BG$}
  {The singularity category of C*BG}}

\begin{definition}
The singularity category $\Dsg(C^*BG)$ is defined to be the Verdier
quotient $\Db(C^*BG)/\Thick(C^*BG)$. This is a triangulated category
with an action of $H^*BG$.
\end{definition}

\begin{proposition}\label{pr:Fred}
If $\Omega BG\phat$ has finite dimensional cohomology then its
connected components are contractible.
\end{proposition}
\begin{proof}
We have $\pi_0(\Omega BG\phat)=\pi_1(BG\phat)=G/O^p(G)$. This group
regularly permutes the connected components of $\Omega BG\phat$, and
the connected component of the identity is $\Omega BO^p(G)\phat$. So
we may assume without loss of generality that $G=O^p(G)$, so that
\[ \pi_0\Omega BG\phat=\pi_1 BG\phat=G/O^pG=1. \]
Now choose an embedding
$G\to SU(n)$, and consider the fibre sequence
\[ \Omega BG\phat \to SU(n)\phat \to (SU(n)/G)\phat \to BG\phat \to
  BSU(n)\phat. \]
We have $\pi_1(SU(n)/G)\phat=1$, and so in the spectral sequence
\[ H^*((SU(n)/G)\phat,H^*(\Omega BG\phat)) \Rightarrow
  H^*(SU(n)\phat), \]
the $E_2$ page is a tensor product.
If $\Omega BG\phat$ has finite
dimensional cohomology, and the top class is in degree greater than
zero,  then the top right hand class in the $E_2$ page of the spectral
sequence is a permanent cycle, contradicting the fact that the top
classes of $SU(n)\phat$ and $(SU(n)/G)\phat$ have the same degree.
Thus $\Omega BG\phat$ is acyclic and has trivial fundamental group, so
it is contractible.
\end{proof}

The following theorem is a special case of the discussion in
Example~3.4 of~\cite{Greenlees/Stevenson:2020a}. We include a proof
for the convenience of the reader.

\begin{theorem}\label{th:Dsg=0}
The singularity category $\Dsg(C^*BG)$ is zero if and only if $G$ is
$p$-nilpotent. 
\end{theorem}
\begin{proof}
If $G$ is $p$-nilpotent then $C^*BG$ is quasi-isomorphic to
$C^*B(G/O^pG)$, and so using Proposition~\ref{pr:DbC*BG} and
Theorem~\ref{th:Kinjb-pgroup} we have
\[ \Db(\kG/O^pG)=\KInjb(\kG/O^pG)\simeq\Db(C^*B(G/O^pG))\simeq\Db(C^*BG). \]
This equivalence sends the generator $\ik$
to $C^*BG$, and so we conclude that $\Thick(C^*BG)=\Db(C^*BG)$ and
$\Dsg(C^*BG)=0$.

Conversely, if $\Dsg(C^*BG)=0$ then $C^*BG$ generates
$\Db(C^*BG)$. By Greenlees and Stevenson
\cite[Example~3.4]{Greenlees/Stevenson:2020a}, this means that $C^*BG$
is regular, and $H_*\Omega BG\phat$ is finite dimensional.
Applying Proposition~\ref{pr:Fred}, it follows that the connected
components of $\Omega BG\phat$ are contractible. Thus $BG\phat$ is a
$K(G/O^pG,1)$, and $G$ is $p$-nilpotent.
\end{proof}

\begin{remark}
  Although $\Db(C^*BG)$ is a tensor triangulated category,
  $\Thick(C^*BG)$ is not an ideal, and the quotient $\Dsg(C^*BG)$ is
  not a tensor triangulated category.
\end{remark}

\begin{conjecture}\label{conj:suppDsg}
  The support of $\Dsg(C^*BG)$ is equal to the nucleus $\Theta_G$.
\end{conjecture}

\begin{theorem}\label{th:C*BS=>suppDsg}
Conjecture~\ref{conj:C*BS} for a finite group $G$ implies that the support of $\Dsg(C^*BG)$
is contained in $\Theta_G$.
\end{theorem}
\begin{proof}
This follows from Theorem~\ref{th:nucleus}.
\end{proof}

\begin{theorem}\label{th:triv-nucleus}
If the nucleus $\Theta_G$ is equal to $\varnothing$ or $\{0\}$ then it
is equal to the support of $\Dsg(C^*BG)$, so that
Conjecture~\ref{conj:suppDsg} is true in this case.
\end{theorem}
\begin{proof}
If $\Theta_G\subseteq\{0\}$ then by Corollary~\ref{co:zero-nucleus},
Conjecture~\ref{conj:C*BS} holds. Then by
Theorem~\ref{th:C*BS=>suppDsg}, it follows that the support of
$\Dsg(C^*BG)$ is contained in $\{0\}$. Now applying
Theorem~\ref{th:Dsg=0}, it is equal to $\varnothing$ is $G$ is
$p$-nilpotent, namely if $\Theta_G=\varnothing$, and it is equal to
$\{0\}$ otherwise.
\end{proof}

\begin{remark}
  Equality of the support of $\Dsg(C^*BG)$ with the nucleus, given
  Conjecture~\ref{conj:C*BS}, seems more difficult to prove if the
  nucleus is non-trivial, but at
  least the methods developed in~\cite{Benson:1995a} enable us to
  reduce to the case of a point in the nucleus which has $G$ as its
  diagonaliser, in the terminology of that paper.
\end{remark}

\bibliographystyle{amsplain}
\bibliography{../repcoh}

\end{document}